\documentclass[twoside,12pt]{article}
\usepackage[]{amsmath,amsthm,amssymb}
\usepackage[latin1]{inputenc}

\begin{document}

\title{{\Large\bf On the orthogonality and convolution orthogonality via the Kontorovich-Lebedev  transform}}

\author{Semyon  YAKUBOVICH}
\maketitle

\markboth{\rm \centerline{ Semyon   YAKUBOVICH}}{}
\markright{\rm \centerline{ Orthogonality via the Kontorovich-Lebedev transform }}

\begin{abstract} {Notions of the orthogonality and convolution orthogonality are explored with the use of the Kontorovich-Lebedev transform and its convolution.
New classes of the corresponding orthogonal polynomials and functions are investigated.  Integral representations, orthogonality relations  and explicit expressions are established.}

\end{abstract}
\vspace{4mm}



\vspace{4mm}

\section {Introduction and preliminary results}

Let $f$ be a complex-valued function defined on $\mathbb{R}_+ \equiv (0,\infty).$ The Kontorovich-Lebedev transform [12], [13] is defined by the following integral

$$(Ff)(\tau)= \int_0^\infty K_{i\tau}(x) f(x) dx,\quad \tau \in \mathbb{R}_+.\eqno(1.1)$$
Here $K_{i\tau}(x)$ is  the modified Bessel function of the second kind,  or the Macdonald function of the argument $x >0$ and the pure
imaginary subscript $i\tau$ (see [5], Vol. II).  It can be defined by the Fourier cosine 
transform
$$K_{i\tau}(x)= \int_{0}^\infty e^{-x\cosh u}\cos(\tau u)\  du,\ x >0, \eqno(1.2)$$
and,  reciprocally,  by the inversion formula we immediately obtain
$$e^{-x\cosh u}= {2\over \pi}\int_{0}^\infty K_{i\tau}(x)\cos(\tau u)\ d\tau.\eqno(1.3)$$
Moreover, it is an eigenfunction for the differential operator
$$ {\mathcal A}\equiv x^2-  x{d \over dx} x {d \over dx},\eqno(1.4)$$
i.e. we have
$${\mathcal A} K_{i\tau}(x) = \tau^2 \ K_{i\tau}(x).\eqno(1.5)$$
The modified Bessel function has  the asymptotic behavior with
respect to $x$ [5], Vol. II 
$$ K_{\nu}(x) = \left( \frac{\pi}{2x} \right)^{1/2} e^{-x} [1+
O(1/x)], \qquad x \to +\infty,\eqno(1.6)$$
$$ K_{\nu}(x) = O\left ( x^{-|{\rm Re}\nu|}\right), \ x \to 0,\eqno(1.7)$$
$$K_0(x) = O(-\log x), \ x \to 0 \eqno(1.8)$$
and with respect to the index $\nu= i\tau$
$$K_{i\tau}(x)= O\left({e^{-\pi\tau/2}\over \sqrt\tau}\right), \
\tau \to +\infty.\eqno(1.9)$$
The following uniform inequality will be useful in the sequel (see [12], formula (1.100))

$$\left| K_{i\tau}(x) \right| \le e^{-\delta\tau} K_0\left(x\cos(\delta)\right),\quad \delta \in \left[ 0, {\pi\over 2} \right].\eqno(1.10)$$

As is known [13], the Kontorovich-Lebedev transform (1.1) extends to a bounded invertible isometric map

$$F: L_2\left(\mathbb{R}_+; x dx\right) \to  L_2\left(\mathbb{R}_+;  {2\over \pi^2} \tau \sinh(\pi\tau) d\tau\right)$$
and the inversion formula holds

$$xf(x)= {2\over \pi^2} \int_0^\infty \tau\sinh(\pi\tau) K_{i\tau}(x) (Ff)(\tau) d\tau,\quad x >0,\eqno(1.11)$$
where integrals (1.1), (1.11) converge with respect to norms of the image spaces.  Moreover, the Parseval equality takes place

$$\int_0^\infty  \left|f(x)\right|^2 x dx = {2\over \pi^2} \int_0^\infty  \tau \sinh(\pi\tau) \left| (Ff)(\tau)\right|^2 d\tau\eqno(1.12)$$
and, via the parallelogram identity,  the generalized Parseval equality

$$ \int_0^\infty  f(x) \overline{g(x)} x dx = {2\over \pi^2} \int_0^\infty  \tau \sinh(\pi\tau) (Ff)(\tau) \overline{(Fg)(\tau)} d\tau,\eqno(1.13)$$
where $Fg$ is the Kontorovich-Lebedev transform (1.1) of a function $g \in  L_2\left(\mathbb{R}_+; x dx\right)$.  

According to [11], [12], the convolution $f *g$ of two functions $f, g$ from the space $L_1\left(\mathbb{R}_+; K_0( p x) dx\right),\ 0 < p \le 1$ related to the Kontorovich-Lebedev transform is given by the formula

$$(f*g)(x)= {1\over 2x} \int_0^\infty\int_0^\infty e^{- {y^2+t^2\over 2yt} x - {yt\over 2x}} f(y)g(t) dydt,\quad x >0.\eqno(1.14)$$
We have

{\bf Theorem 1.} {\it Let $f, g \in  L_1\left(\mathbb{R}_+; K_0( p x) dx\right),\ 0 < p \le 1$. Then the convolution $f*g \in  L_1\left(\mathbb{R}_+; K_0( p^2 x) dx\right) $  and satisfies the Young-type inequality

$$||f*g||_{L_1\left(\mathbb{R}_+; K_0( p^2 x) dx\right)} \le  || f ||_{L_1\left(\mathbb{R}_+; K_0( p x) dx\right)} ||g||_{L_1\left(\mathbb{R}_+; K_0( p x) dx\right)}.\eqno(1.15)$$
Moreover, the Kontorovich-Lebedev transform $(1.1)$ of the convolution $(1.14)$ is the product of the Kontorovich-Lebedev transforms, i.e.

$$F (f*g)(\tau) = (Ff)(\tau) (Fg)(\tau),\eqno(1.16)$$
and when $p \in (0, 1/2)$ the Parseval-type equality holds for all $ x>0$}

$$(f*g)(x) =  {2\over  x \pi^2} \int_0^\infty \tau\sinh(\pi\tau) K_{i\tau}(x) (Ff)(\tau) (Fg)(\tau) d\tau.\eqno(1.17)$$

\begin{proof} In fact, the existence of the convolution (1.14) as a function of the space  $ L_1\left(\mathbb{R}_+;\right.$ \\  $\left. K_0( p^2 x) dx\right) $ and inequality (1.15) follows  from the following estimate

$$ ||f*g||_{L_1\left(\mathbb{R}_+; K_0( p^2 x) dx\right)} =    \int_0^\infty {K_0( p^2 x) \over 2x} \left| \int_0^\infty\int_0^\infty e^{- {y^2+t^2\over 2yt} x - {yt\over 2x}} f(y)g(t) dydt \right| dx $$

$$ \le   \int_0^\infty {K_0( p^2 x) \over 2x}  \int_0^\infty\int_0^\infty e^{- {y^2+t^2\over 2yt} x - {yt\over 2x}} \left| f(y)g(t)\right| dydt dx $$

$$= \int_0^\infty {K_0(  x) \over 2x}  \int_0^\infty\int_0^\infty e^{- {y^2+t^2\over 2yt p^2} x - {yt p^2\over 2x}} \left| f(y)g(t)\right| dydt dx $$

$$\le \int_0^\infty {K_0(  x) \over 2x}  \int_0^\infty\int_0^\infty e^{- {y^2+t^2\over 2yt } x - {yt p^2\over 2x}} \left| f(y)g(t)\right| dydt dx $$

$$=  \int_0^\infty K_0( p y)  |f(y)| dy \int_0^\infty K_0( p t)  |g(t)| dt  =   || f ||_{L_1\left(\mathbb{R}_+; K_0( p x) dx\right)} ||g||_{L_1\left(\mathbb{R}_+; K_0( p x) dx\right)},$$
where the interchange of the order of integration is permitted via Fubini's theorem and the integral with respect to $x$  is calculated by the Macdonald product formula for the modified Bessel functions (see [13], formula (1.103)).  In order to prove formula (1.16) we apply the Kontorovich-Lebedev transform to the convolution, change the order of integration by Fubini's theorem and use the Macdonald formula, taking into account the embeddings 

$$ L_1\left(\mathbb{R}_+; K_0( p^2 x) dx\right) \subseteq L_1\left(\mathbb{R}_+; K_0( p x) dx\right) \subseteq  L_1\left(\mathbb{R}_+; K_0( x) dx\right).$$ 
Finally, we establish representation (1.17) for  convolution (1.14). To do this, we appeal to the  following index integral for the convolution kernel (see [13], formula (4.36))

$$ \exp \left(- {1\over 2} \left( {y^2+t^2\over yt} x + {yt\over x}\right) \right) = {4\over \pi^2} \int_0^\infty  \tau\sinh(\pi\tau) K_{i\tau}(x)  K_{i\tau}(y) K_{i\tau}(t) d\tau.\eqno(1.18)$$ 
Hence, employing inequality (1.10), we deduce

$$  \exp \left(- {1\over 2} \left( {y^2+t^2\over yt} x + {yt\over x}\right) \right) \le   {4\over \pi^2} \ K_{0}(x\cos(\delta)) K_{0}(y\cos(\delta))K_{0}(t\cos(\delta))  $$

$$\times  \int_0^\infty  \tau\sinh(\pi\tau) e^{-3\delta \tau} d\tau =  {24 \delta \over \pi (9\delta^2- \pi^2)} \ K_{0}(x\cos(\delta)) K_{0}(y\cos(\delta))K_{0}(t\cos(\delta)),$$
which gives the inequality for all $(x,y,t) \in \mathbb{R}^3_+,\  \delta \in \left({\pi\over 3},\ {\pi\over 2} \right)$

$$  \exp \left(- {1\over 2} \left( {y^2+t^2\over yt} x + {yt\over x}\right) \right) \le   {24  \delta \over \pi (9\delta^2- \pi^2)} \ K_{0}(x\cos(\delta)) K_{0}(y\cos(\delta))K_{0}(t\cos(\delta)).\eqno(1.19)$$
Hence, plugging  the right-hand side of (1.18) in (1.14)  and changing the order of integration owing to inequality (1.19) and conditions of the theorem for $p= \cos(\delta) \in (0, 1/2)$,
we get representation (1.17) of the convolution $f*g$.

\end{proof}

{\bf Lemma 1}. {\it Let $f  \in  L_1\left(\mathbb{R}_+; K_0( p x)  x^{-1} dx\right),\  ,\ 0 < p \le 1$.  Then the Kontorovich-Lebedev transform $(1.1)$ $(F g)(\tau)$ of the function $g(x)= f(x)/x$ is the composition of the Fourier cosine and Laplace  transforms, i.e.

$$(Fg)(\tau)= \int_0^\infty \cos(\tau u) \int_0^\infty  e^{-x\cosh u} f(x) { dx du\over x},\quad \tau >0.\eqno(1.20)$$
Moreover, if, in addition,  $f$  belongs to the space  $L_r\left(\mathbb{R}_+; K_0( p x)  x^{-a} dx\right),\\ r, a >1,\ p \in (0, 1)$, then it can be written in the form

$$(Fg)(\tau)= {1\over \tau} \int_0^\infty \sin(\tau u) \sinh(u) \int_0^\infty  e^{-x\cosh u} f(x) dx du\eqno(1.21)$$
and}
$$\lim_{u\to \infty} \sinh(u) \int_0^\infty  e^{-x\cosh u} f(x) dx = 0.\eqno(1.22)$$

\begin{proof} Since  $f  \in  L_1\left(\mathbb{R}_+; K_0( p x)  x^{-1} dx\right)$  we have the estimate

$$ \int_0^\infty \left| \cos(\tau u) \int_0^\infty  e^{-x\cosh u} f(x) { dx\over x}\right| \le  \int_0^\infty du  \int_0^\infty  e^{-x\cosh u} |f(x)| { dx \over x}$$

$$=  \int_0^\infty  K_0(x) |f(x)| { dx \over x} \le ||f||_{ L_1\left(\mathbb{R}_+; K_0( p x)  x^{-1} dx\right) } < \infty.$$
Therefore formula (1.20) follows immediately via Fubini's theorem and integral representation (1.2).  Then integrating by parts in the integral with respect to $u$ and eliminating the integrated terms due to the absolute and uniform convergence by $u \in \mathbb{R}_+$ of the integral with respect to $x$, we obtain (1.21), differentiating by $u$ under the integral sign. This is, indeed,  allowed due to the estimate for some $N >0$ large  enough

$$ \sinh(u) \int_N^\infty  e^{-x\cosh u} |f(x)| dx =    \sinh(u) \int_N^\infty  e^{-(1-p) x\cosh u}  e^{- p x\cosh u} |f(x)| dx $$

$$\le  \sinh(u) \left(\int_N^\infty  e^{-(1-p)r^\prime x\cosh u}  x^{a r^\prime/r} dx\right)^{1/r^\prime}  \left( \int_N^\infty  e^{- p r x\cosh u} {|f(x)|^r\over x^a} dx \right)^{1/r}$$

$$\le   \Gamma^{1/r^\prime} \left( {a r^\prime\over r} +1\right)  { \sinh(u) \over [ (1-p) r^\prime \cosh(u)]^{ 1+ (a-1)/r} }  \left( \int_N^\infty  e^{- p r x\cosh u} {|f(x)|^r\over x^a} dx \right)^{1/r}$$

$$\le   \Gamma^{1/r^\prime} \left( {a r^\prime\over r} +1\right)  \left[ (1-p) r^\prime \right]^{ -1- (a-1)/r}  $$

$$\times \sup_{x \ge N} \left[ { e^{- p r x}\over K_0(px)}\right]^{1/r} \left( \int_N^\infty  K_0(px)  {|f(x)|^r\over x^a} dx \right)^{1/r} \to 0,\quad N \to \infty,$$
where $r^\prime= r/(r-1)$ and $\Gamma(z)$ is the Euler gamma-function [5], Vol. I.   Finally, condition (1.22) is an immediate consequence of the inequality

$$ \sinh(u) \int_0^\infty  e^{-x\cosh u} |f(x)| dx \le   \Gamma^{1/r^\prime} \left( {a r^\prime\over r} +1\right)  { [\cosh(u)]^{(1-a)/r}  \over [ (1-p) r^\prime ]^{ 1+ (a-1)/r} } $$
$$\times \sup_{x \ge N} \left[ { e^{- p r x}\over K_0(px)}\right]^{1/r}  ||f||_{ L_r\left(\mathbb{R}_+; K_0( p x)  x^{-a} dx\right) }.$$

\end{proof}

{\bf Corollary 1}. {\it Let  $f, g \in  L_1\left(\mathbb{R}_+; K_0( p x) dx\right),\ 0 < p < 1/2$ and $\omega \in L_1\left(\mathbb{R}_+; K_0( p x)  x^{-1} dx\right).$ Then the following equality takes place

$$ \int_0^\infty (f*g)(x) \omega(x) dx =  {2\over   \pi^2} \int_0^\infty \tau\sinh(\pi\tau) (Ff)(\tau) (Fg)(\tau) q(\tau) d\tau,\eqno(1.23)$$
where}

$$q(\tau)=  \int_0^\infty K_{i\tau}(x) \omega(x) {dx\over x}.\eqno(1.24)$$

\begin{proof}  The proof follows immediately, multiplying both sides of (1.17) by $\omega$ and integrating over $\mathbb{R}_+$. The interchange of the order of integration on
 the right-hand side of the obtained equality is allowed via Fubini's theorem and  inequality (1.19).

\end{proof}

Let $\omega$ be a positive function such that $\omega \in L_1\left( (0,1);  x^{-2} dx\right) \cap L_1\left( (1,\infty);  x^{-1} dx\right)$.  Making use  integral representations for the modified Bessel function  via the integration by parts [5], Vol. II

$$K_{i\tau}(x)= {1\over \sinh(\pi\tau/2)} \int_0^\infty \sin(\tau u) \sin(x\sinh(u)) du $$

$$=   {1\over  x \sinh(\pi\tau/2)} \int_0^\infty  {\cos(x\sinh(u)) \over \cosh^2(u)} \left[ \tau \cos(\tau u) \cosh(u) - \sin(\tau u) \sinh(u) \right] du,\eqno(1.25)$$
we substitute the first  integral in (1.25) into (1.24), having the equality

$$ q(\tau)=  {1\over \sinh(\pi\tau/2)}  \int_0^\infty   \int_0^\infty \sin(\tau u) \sin(x\sinh(u))  \omega(x) {dx du\over x}.\eqno(1.26)$$
Our goal is to justify the interchange of the order of integration in (1.26), proving the formula

$$ q(\tau)=  {1\over \sinh(\pi\tau/2)}   \int_0^\infty \sin(\tau u) du  \int_0^\infty  \sin(x\sinh(u))  \omega(x) {dx\over x}.\eqno(1.27)$$
Indeed, we write the integral in (1.27) in the form

$$   \lim_{N\to \infty}  \int_0^N  \sin(\tau u) du  \int_0^\infty  \sin(x\sinh(u))  \omega(x) {dx\over x} $$

$$=   \lim_{N\to \infty}    \int_0^\infty  \int_0^N \sin(\tau u) \sin(x\sinh(u))   \omega(x) { du dx\over x},\eqno(1.28)$$
where the interchange of the order of integration is due to the dominated convergence theorem by virtue of the estimate 

$$  \int_0^N \left| \sin(\tau u) \right| du  \int_0^\infty \left| \sin(x\sinh(u)) \right|  \omega(x) {dx\over x} \le N \int_0^\infty  \omega(x) {dx\over x} < \infty.$$ 
Now,   in order to pass to the limit under the integral sign on the right-hand side  of (1.28), we appeal to the same justification  as above, employing the following estimate via the integration by parts (cf. (1.25))

$$ \left| \int_0^N \sin(\tau u) \sin(x\sinh(u))  du \right| = {1\over x} \left|-  {\sin(N\tau) \cos(x \sinh(N)) \over \cosh(N)} \right.$$

$$\left. + \int_0^N  {\cos(x\sinh(u)) \over \cosh^2(u)} \left[ \tau \cos(\tau u) \cosh(u) - \sin(\tau u) \sinh(u) \right] du\right| $$

$$\le  {1\over x} \left[ 1+   \int_0^\infty  {1 \over \cosh^2(u)} \left[ \tau \cosh(u) + \sinh(u) \right] du \right]$$
and the condition $\omega \in L_1\left(\mathbb{R}_+;  x^{-2} dx\right).$ This proves (1.27).  Further, appealing to Theorem 123 in [10] about the non-negativeness of the Fourier sine transform, 
we are ready, as a direct consequence,   to formulate the following result.

{\bf Lemma 2}. {\it Let  $\omega$ be a positive non-increasing function over $\mathbb{R}_+$ such that $\omega(x)= o(x),\ x \to \infty, \omega \in L_1\left( (0,1);  x^{-2} dx\right) \cap L_1\left( (1,\infty);  x^{-1} dx\right)$ and let the function

$$\varphi(u)=  \int_0^\infty  \sin(x u)  \omega(x) {dx\over x}\eqno(1.29)$$
is non-increasing over $\mathbb{R}_+$. Then the weight function $(1.24)$ $q(\tau) \ge 0,\ \tau \in \mathbb{R}_+$.}

\begin{proof}    It is easily seen from the conditions of the lemma and Theorem 123 in [10] that $\varphi(u)$ as the Fourier sine transform (1.29) is non-negative. Then since $\varphi(\sinh(u))$ is non-increasing over $\mathbb{R}_+$, integrable over $(0,1)$ owing to the estimate

$$\int_0^1 \left| \varphi(\sinh(u)) \right| du \le ||\omega||_{ L_1\left(\mathbb{R}_+;  x^{-1} dx\right)}$$
and tends to zero at infinity due to the Riemann-Lebesgue lemma, its Fourier sine transform (1.27)  is non-negative as well  via Theorem 123 in [10]. Hence $q(\tau) \ge 0,\ \tau \in \mathbb{R}_+$. 
\end{proof} 

Returning to  equality (1.23), we observe that for positive $\omega$ and non-negative $q$ its  left-hand side satisfies properties of the inner product and forms the so-called convolution Hilbert spaces studied in [13, Chapter 4]. This means that one can consider the corresponding convolution orthogonality of functions from the associated Lebesgue spaces. The notion of the convolution orthogonality was introduced for the first time in [2] for the Laplace convolution, and the discrete case was investigated in [1]. Our goal is to explore in the sequel concrete orthogonal sequences  of functions and polynomials with respect to the convolution (1.14) and Parseval-type equalities (1.13), (1.23). To do this we will employ, in particular, Wilson's and Continuous Dual Hahn polynomials [11]. We note that in [6], [9] similar problems were examined, involving  Fourier and Fourier-Jacobi transforms.   

Finally, in this section, let $\alpha \in \mathbb{R}_+,\ \beta \in \mathbb{R} \backslash\{0\}$ and  consider the modified Kontorovich-Lebedev transform by the formula

$$(F_{\alpha,\beta} f)(\tau)=  \int_0^\infty K_{i\tau}(\alpha x^\beta) f(x) dx.\eqno(1.30)$$
We have (see (1.1)) $(F_{1,1} f)(\tau) \equiv (Ff)(\tau)$.  Furthermore, in the same manner as above, changing variables and functions, one can show that the modified Kontorovich-Lebedev transform  (1.30) extends to a bounded invertible isometric map
$$F_{\alpha,\beta}: L_2\left(\mathbb{R}_+; x dx\right) \to  L_2\left(\mathbb{R}_+;  {2 |\beta| \over \pi^2} \tau \sinh(\pi\tau) d\tau\right),\eqno(1.31)$$
 the inversion formula holds

$$xf(x)= {2 |\beta| \over \pi^2} \int_0^\infty \tau\sinh(\pi\tau) K_{i\tau}(\alpha x^\beta) (F_{\alpha,\beta}f)(\tau) d\tau,\quad x >0,\eqno(1.32)$$
where integrals (1.30), (1.31) converge with respect to norms of the image spaces, and the Parseval equality is valid 

$$ \int_0^\infty  f(x) \overline{g(x)} x dx = {2 |\beta| \over \pi^2} \int_0^\infty  \tau \sinh(\pi\tau) (F_{\alpha,\beta} f)(\tau) \overline{(F_{\alpha,\beta}g)(\tau)} d\tau.\eqno(1.33)$$

\section{The use of the Wilson  and Continuous Dual Hahn polynomials}

With a slight modification of the Kontorovich-Lebedev transform (1.1) it becomes an automorphism on the vector space of polynomials $\mathcal{P}$ (see [7], [8]).   Moreover, it involves orthogonal and multiple orthogonal polynomials or $d$-orthogonal ones. For instance, the Continuous Dual Hahn polynomials appear as the Kontorovich-Lebedev transform of a 2-orthogonal sequence of Laguerre type.  Our goal here is to extend Parseval-type equalities (1.13), (1.23) to  different orthogonal polynomials  and functions, generating new  orthogonal and convolution orthogonal systems related to the Kontorovich-Lebedev transform and its modifications.  It is straightforward to verify,  appealing to the asymptotic behavior (1.6), (1.7), (1.8) of the modified Bessel function, that $\mathcal{P} \subset   L_1\left(\mathbb{R}_+; K_0( p x) dx\right),\ 0 < p \le 1.$ This means, that the Parseval-type equality (1.23) holds for polynomial functions $f,g$. 
However, the Parseval equality (1.13) to be valid for a polynomial $f$ should be reestablished  under some sufficient conditions on a function $g$. To do this, we shall prove 

{\bf Theorem 2.} {\it Let $f \in \mathcal{P}$ and  $g$ be such that its Kontorovich-Lebedev transform $Fg  \in L_2\left(\mathbb{R}_+; \tau\sinh((\pi+\mu)\tau) d\tau \right)$, for some $\mu >0$.   Then equality $(1.13)$ holds.}

\begin{proof}  Since evidently,  $L_2\left(\mathbb{R}_+; \tau\sinh((\pi+\mu)\tau) d\tau \right)  \subset   L_2\left(\mathbb{R}_+; \tau\sinh(\pi\tau) d\tau \right),\ \mu >0$,  we get from (1.12) that $g \in L_2\left(\mathbb{R}_+;  x dx \right)$ and  integral (1.11) for $g$

$$xg(x)= {2\over \pi^2} \int_0^\infty \tau\sinh(\pi\tau) K_{i\tau}(x) (Fg)(\tau) d\tau,\eqno(2.1)$$
converges in the mean square sense with respect to the norm in  $L_2\left(\mathbb{R}_+;  x dx \right)$. But the estimate (see (1.10)) 

$$ \int_0^\infty \tau\sinh(\pi\tau) \left| K_{i\tau}(x) (Fg)(\tau) \right| d\tau \le K_0(x\cos(\delta))  \int_0^\infty \tau\sinh(\pi\tau) e^{-\delta\tau} \left| (Fg)(\tau) \right| d\tau$$

$$\le K_0(x\cos(\delta))  \left( \int_0^\infty {\tau\sinh^2 (\pi\tau) e^{-2\delta \tau}\over  \sinh((\pi+\mu)\tau)} \ d\tau \right)^{1/2} || Fg||_{L_2\left(\mathbb{R}_+; \tau\sinh((\pi+\mu)\tau) d\tau \right)}$$
$$=   C_{\delta,\mu} K_0(x\cos(\delta))  || Fg||_{L_2\left(\mathbb{R}_+; \tau\sinh((\pi+\mu)\tau) d\tau \right)},$$
where $C_{\delta,\mu} >0$ is the constant

$$ C_{\delta,\mu} =   \left( \int_0^\infty {\tau\sinh^2 (\pi\tau) e^{-2\delta \tau}\over  \sinh((\pi+\mu)\tau)} \ d\tau \right)^{1/2},\quad \delta \in \left( \max\left( 0, {\pi -\mu\over 2}\right), {\pi\over 2} \right),$$
guarantees the existence of (2.1) as  a Lebesgue integral for all $ x>0$.  Therefore, substituting it in the left-hand side of (1.13), we change the order of integration by Fubini's theorem owing to the estimate

$$  \int_0^\infty  |f(x)|  \int_0^\infty \tau\sinh(\pi\tau) \left|K_{i\tau}(x) (Fg)(\tau) \right| d\tau dx $$

$$\le  C_{\delta,\mu}   || Fg||_{L_2\left(\mathbb{R}_+; \tau\sinh((\pi+\mu)\tau) d\tau \right)} \int_0^\infty  |f(x)| K_0(x\cos(\delta)) dx  < \infty,$$
to complete the proof of Theorem 2. 

\end{proof}

We will show below that the Parseval equality (1.13) generates various systems of orthogonal polynomials and functions, and it is closely related, in particular,  to the Wilson and Continuous Dual Hahn orthogonalities.  Precisely, following  [11], we define the Wilson real orthogonality by the equality

$$\int_0^\infty \left| \frac{ \Gamma(a+it)\Gamma(b+it)\Gamma(c+it)\Gamma(d+it)}{\Gamma(2it)}\right|^2 W_n( t^2) W_m(t^2) dt $$

$$= 2\pi \delta_{n,m} \frac{ n! \Gamma(n+a+b)\Gamma(n+a+c)\Gamma(n+a+d)\Gamma(n+b+c)\Gamma(n+b+d)\Gamma(n+c+d)}
{\Gamma(2n+a+b+c+d) (n+a+b+c+d-1)_n},$$
where $\delta_{n,m}$ is the Kronecker  delta, $(z)_n,\ n \in \mathbb{N}_0$ is the Pochhammer symbol, $a,b,c,d >0,\ i$ is the imaginary unit and 

$$W_n(t^2)\equiv  W_n(t^2; a,b,c,d) = (a+b)_n(a+c)_n (a+d)_n$$

$$\times  {}_4F_3 \left( -n, n+a+b+c+d-1, a+it, a-it; a+b, a+c, a+d; 1\right)\eqno(2.2)$$
are Wilson's polynomials being expressed in terms of the generalized hypergeometric function ${}_4F_3$ (cf. [5], Vol. I).  It is shown that $W_n$ is symmetric in all four parameters $a,b,c,d$.

The Continuous Dual Hahn orthogonality is defined accordingly [11]

$$\int_0^\infty \left| \frac{ \Gamma(a+it)\Gamma(b+it)\Gamma(c+it)}{\Gamma(2it)}\right|^2 S_n( t^2) S_m(t^2) dt $$

$$= 2\pi \delta_{n,m} \  n! \Gamma(n+a+b)\Gamma(n+a+c)\Gamma(n+b+c),$$
where 

$$S_n(t^2)\equiv  S_n(t^2; a,b,c) = (a+b)_n(a+c)_n $$

$$\times  {}_3F_2 \left( -n, a+it, a-it; a+b, a+c; 1\right)\eqno(2.3)$$
are Continuous Dual Hahn  polynomials in terms of the ${}_3F_2$ generalized hypergeometric function.  These polynomials are symmetric in all three parameters. 

Let us consider the classical associated Laguerre orthogonality [5], Vol. II of the corresponding  polynomials  $ L_n^\alpha(x)$

$$\int_0^\infty x^\alpha e^{-x} L_n^\alpha(x) L_m^\alpha(x) dx = {\Gamma(n+\alpha+1)\over n!} \delta_{n,m},\quad \alpha > -1.\eqno(2.4)$$
We will see as this orthogonality generates other orthogonalities of polynomials and functions, making use the Parseval equality  (1.13)  and mapping properties of the Kontorovich-Lebedev transform.  In fact,  denoting the following Kontorovich-Lebedev integral by

$$ F_n\left(\tau; \alpha,\beta,\mu, \eta\right)= \int_0^\infty x^\beta e^{- \mu x} L_n^\alpha(x) K_{i\tau} (\eta x) dx,\quad \alpha, \beta > -1,\ \mu \ge 0,\ \eta >0,\eqno(2.5)$$
we let $\alpha=\beta+\gamma+1,\ \beta, \gamma > -1$  to satisfy mapping $L_2$-properties of the Kontorovich-Lebedev transform (1.1) (see (1.12), (1.33) and Theorem 2) in order to apply  the Parseval equality (1.13) to the left-hand side of (2.4). As the result we obtain ($ 0 < \mu \le 1$)

$$  \int_0^\infty  \tau \sinh(\pi\tau) F_n\left(\tau; \alpha,\beta,\mu,\mu\right) F_m\left(\tau; \alpha,\gamma, 1-\mu,\mu\right) d\tau = {\pi^2\over 2 n!} \ \Gamma(n+\alpha+1) \delta_{n,m}.\eqno(2.6)$$
But the function $F_n\left(\tau; \alpha,\beta,\mu,\mu\right)$ can be calculated explicitly, employing Entry 3.14.3.1 in [3]

$$\int_0^\infty  x^{s-1} e^{-\mu x} K_{i\tau} (\mu x) dx  = {\sqrt\pi\over (2\mu)^s} \frac{\Gamma(s-i\tau)\Gamma(s+i\tau)}{ \Gamma(s+1/2)},\ {\rm Re} s >0\eqno(2.7)$$
and the closed form expression of the associated Laguerre polynomials

$$L_n^\alpha(x) = \sum_{k=0}^n (-1)^k \binom{n+\alpha}{n-k} {x^k\over k!}.\eqno(2.8)$$
Hence we derive from (2.4)

$$ F_n\left(\tau; \alpha,\beta,\mu,\mu\right) = {\sqrt\pi\over (2\mu)^{\beta+1}} \sum_{k=0}^n {(-1)^k \over k!} \binom{n+\alpha}{n-k} (2\mu)^{-k} \frac{\Gamma(\beta+k+1-i\tau)\Gamma(\beta+k+1+i\tau)}{ \Gamma(\beta+k+3/2)}$$

$$=  {\sqrt\pi\over (2\mu)^{\beta+1} n!} {\left| \Gamma(\beta+1-i\tau)\right|^2 (\alpha+1)_n \over \Gamma(\beta+3/2)} $$

$$\times {}_3F_2 \left( - n,\  \beta+1+i\tau,\  \beta+1-i\tau;\  \alpha+1,\ \beta+ {3\over 2};\   {1\over 2\mu} \right),\eqno(2.9)$$
where the latter ${}_3F_2$ hypergeometric function is a polynomial in $\tau$, slightly generalizing  the Continuous Dual Hahn polynomial (2.3).  In particular, letting $\mu=1/2$ in (2.6) and (2.9), employing the reflection and duplication formulae for gamma-function [5], Vol. I, we end up with a particular case of the Continuous Dual Hahn orthogonality, namely,

$$    \int_0^\infty \left| \frac{ \Gamma(\beta+1+it)\Gamma(\gamma+1+it)}{  \Gamma(it)}\right|^2 S_n \left( t^2; \beta+1, \gamma+1, {1\over 2}\right)  S_m\left( t^2; \beta+1, \gamma+1, {1\over 2}\right) dt $$

$$=  {n! \over 2} \ \delta_{n,m} \  \Gamma(n+\beta+\gamma+2) \Gamma \left(n+\beta+{3\over 2}\right)  \Gamma\left(n+\gamma+{3\over 2}\right).\eqno(2.10)$$
On the other hand, we calculate $F_m\left(\tau; \alpha,\gamma, 1-\mu,\mu\right) $, appealing to (2.8) and Entry 3.14.3.5 in [3], which can be written  in terms of the Gauss hypergeometric function $ {}_2F_1$ [5], Vol. I.  Precisely, we have

$$\int_0^\infty  x^{s-1}  e^{-(1-\mu) x} K_{i\tau} (\mu x) dx = \sqrt{\pi} \  \frac{  \Gamma(s-i\tau)\Gamma(s+i\tau)}{ (2\mu)^s \Gamma(s+1/2)}$$

$$\times  {}_2F_1 \left( s+i\tau,\ s-i\tau; \ s+ {1\over 2}; \ 1- {1\over 2\mu} \right).\eqno(2.11)$$
Then we find, finally,  the following expressions

$$ F_m\left(\tau; \alpha,\gamma, 1-\mu,\mu\right)  =  \int_0^\infty x^\gamma  e^{- (1-\mu) x} L_m^\alpha(x) K_{i\tau} (\mu x) dx $$

$$=   \sqrt \pi \frac{ |\Gamma(\gamma+1+i\tau)|^2\  (1+\alpha)_m}{ (2\mu)^{\gamma+1} \ m! \ \Gamma(\gamma+3/2)}\sum_{k=0}^m  \frac{ (-m)_k  (\gamma+1 -i\tau)_k (\gamma+1+i\tau)_k}{ (2\mu)^k\  k! \  (1+\alpha)_k  (\gamma+3/2)_k}$$

$$\times  {}_2F_1 \left( k+\gamma+1+i\tau,\ k+\gamma+1-i\tau; \ k+\gamma+ {3\over 2}; \ 1- {1\over 2\mu} \right)$$

$$ = \sqrt{ {\pi\over 2\mu}} \  {(1+\alpha)_m\over  m!}\  |\Gamma(\gamma+1+i\tau)|^2 \sum_{k=0}^m  \frac{ (-m)_k  (\gamma+1 -i\tau)_k (\gamma+1+i\tau)_k}{  k! \  (1+\alpha)_k }$$

$$\times  (1-2\mu)^{-(k+\gamma+1/2)/2}  P^{-k-\gamma-1/2}_{i\tau-1/2} \left({1\over \mu}- 1\right),\eqno(2.12)$$
where $P_\nu^\lambda (z)$ is the associated Legendre function [5], Vol. I.   An important special case of the orthogonality (2.6) is $\mu=1$.  In fact,  recalling (2.8) and using Entry  3.14.1.3 in [3], we get

$$ F_{2m}\left(\tau; \alpha,\gamma,  0, 1\right)  =  \int_0^\infty x^\gamma   L_{2m}^\alpha(x) K_{i\tau} (x) dx = {(1+\alpha)_{2m}\over (2m)!} \sum_{k=0}^{2m} { 2^{k+\gamma-1} (-2m)_k\over k! (1+\alpha)_k }  $$

$$\times \Gamma\left({k+\gamma+1+i\tau\over 2}\right)\Gamma\left({k+\gamma+1-i\tau\over 2}\right) =  {(1+\alpha)_{2m}\over (2m)!} \left[ \left|\Gamma\left({\gamma+1+i\tau\over 2}\right) \right|^2\right.$$

$$\times \sum_{k=0}^{m} { 2^{2k+\gamma-1} (-2m)_{2k}\over (2k)! (1+\alpha)_{2k} } \ \left({\gamma+1+i\tau\over 2}\right)_k  \left({\gamma+1-i\tau\over 2}\right)_k$$

$$+ \left. \left|\Gamma\left({\gamma+2+i\tau\over 2}\right) \right|^2\sum_{k=0}^{m-1} { 2^{2k+\gamma} (-2m)_{2k+1}\over (2k+1)! (1+\alpha)_{2k+1} } \ \left({\gamma+2+i\tau\over 2}\right)_k  \left({\gamma+2-i\tau\over 2}\right)_k\right]$$

$$ = {(1+\alpha)_{2m}\  2^\gamma \over (2m)!} \left[  {1\over 2} \left|\Gamma\left({\gamma+1+i\tau\over 2}\right) \right|^2\right.$$

$$\times \  {}_4F_3 \left( -m, \ {1\over 2} - m,\  {\gamma+1+i\tau\over 2},\ {\gamma+1- i\tau\over 2} ; \ {1\over 2},\  {1+\alpha\over 2},\  1+{ \alpha\over 2};\ 1\right)$$

$$-  { 2m\over 1+\alpha}  \left|\Gamma\left({\gamma+2+i\tau\over 2}\right) \right|^2 $$

$$\left.\times  {}_4F_3 \left( 1 -m, \ {1\over 2} - m,\  {\gamma+2+i\tau\over 2},\ {\gamma+2- i\tau\over 2} ; \ {3\over 2},\  1+{\alpha\over 2},\  {3+ \alpha\over 2};\ 1\right)\right],\eqno(2.13)$$

$$ F_{2m+1}\left(\tau; \alpha,\gamma,  0, 1\right)  =  \int_0^\infty x^\gamma   L_{2m+1}^\alpha(x) K_{i\tau} (x) dx = {(1+\alpha)_{2m+1}\over (2m+1)!} \sum_{k=0}^{2m+1} { 2^{k+\gamma-1} (-2m-1)_k\over k! (1+\alpha)_k }  $$

$$\times \Gamma\left({k+\gamma+1+i\tau\over 2}\right)\Gamma\left({k+\gamma+1-i\tau\over 2}\right) =  {(1+\alpha)_{2m+1}\over (2m+1)!} \left[ \left|\Gamma\left({\gamma+1+i\tau\over 2}\right) \right|^2\right.$$

$$\times \sum_{k=0}^{m} { 2^{2k+\gamma-1} (-2m-1)_{2k}\over (2k)! (1+\alpha)_{2k} } \ \left({\gamma+1+i\tau\over 2}\right)_k  \left({\gamma+1-i\tau\over 2}\right)_k$$

$$+ \left. \left|\Gamma\left({\gamma+2+i\tau\over 2}\right) \right|^2\sum_{k=0}^{m} { 2^{2k+\gamma} (-2m-1)_{2k+1}\over (2k+1)! (1+\alpha)_{2k+1} } \ \left({\gamma+2+i\tau\over 2}\right)_k  \left({\gamma+2-i\tau\over 2}\right)_k\right]$$

$$ = {(1+\alpha)_{2m+1}\  2^\gamma \over (2m+1)!} \left[  {1\over 2} \left|\Gamma\left({\gamma+1+i\tau\over 2}\right) \right|^2\right.$$

$$\times \  {}_4F_3 \left( -m, \  - {1\over 2} - m,\  {\gamma+1+i\tau\over 2},\ {\gamma+1- i\tau\over 2} ; \ {1\over 2},\  {1+\alpha\over 2},\  1+{ \alpha\over 2};\ 1\right)$$

$$-  { 2m+1\over 1+\alpha}  \left|\Gamma\left({\gamma+2+i\tau\over 2}\right) \right|^2 $$

$$\left.\times  {}_4F_3 \left(  -m, \ {1\over 2} - m,\  {\gamma+2+i\tau\over 2},\ {\gamma+2- i\tau\over 2} ; \ {3\over 2},\  1+{\alpha\over 2},\  {3+ \alpha\over 2};\ 1\right)\right].\eqno(2.14)$$

Let us treat the associated Laguerre orthogonality with the aid of the modified Kontorovich-Lebedev transform (1.30), letting $\beta=1/2$.  Then equality (2.6) becomes

$$ \int_0^\infty  \tau \sinh(\pi\tau) G_n\left(\tau; \alpha,\beta,\mu,\eta\right) G_m\left(\tau; \alpha,\gamma, 1-\mu,\eta\right) d\tau
=  { \pi^2 \over n!} \  \Gamma(n+\alpha+1) \delta_{n,m},\eqno(2.15)$$
where $\alpha= \beta+\gamma+1,\ \alpha, \beta,\ \gamma > -1,\  0 < \mu \le 1,\ \eta >0$ and 

$$ G_n\left(\tau; \alpha,\beta,\mu,\eta\right) =  \int_0^\infty x^\beta e^{- \mu x} L_n^\alpha(x) K_{i\tau}\left(\eta \sqrt x\right) dx,\eqno(2.16)$$

$$ G_m\left(\tau; \alpha,\gamma,1-\mu,\eta\right) =  \int_0^\infty x^\gamma e^{- (1-\mu) x} L_m^\alpha(x) K_{i\tau}\left(\eta\sqrt x\right) dx.\eqno(2.17)$$
Then (2.8) and Entry 3.14.3.10 in [3]  suggest the equalities 

$$G_n\left(\tau; \alpha,\beta,\mu,\eta\right)  =  \sum_{k=0}^n {(-1)^k \over k!} \binom{n+\alpha}{n-k} 
\int_0^\infty x^{\beta+k}  e^{- \mu x}  K_{i\tau}\left(\eta \sqrt x\right) dx $$

$$= {(1+\alpha)_n\over \eta \ n! } \ \mu^{- (1/2+\beta)}\  e^{\eta^2/ (8\mu)} \left|\Gamma\left( \beta+1+{i\tau\over 2} \right)\right|^2   \sum_{k=0}^n {(-n)_k \over k!} $$
 
$$\times {\left( \beta+1-{i\tau\over 2} \right)_k \left( \beta+1-{i\tau\over 2} \right)_k\over (1+\alpha)_k} \  W_{ - (1/2+\beta+k), i\tau/2} \left({\eta^2\over 4\mu}\right),\eqno(2.18)$$ 
where $W_{\nu,\mu}(z)$ is the Whittaker function [5], Vol. II,

$$  G_m\left(\tau; \alpha,\gamma,1-\mu,\eta\right) =  {(1+\alpha)_m\over \eta\ m!} \ (1-\mu)^{- (1/2+\gamma}\  e^{\eta^2/ (8(1-\mu))} \left|\Gamma\left( \gamma+1-{i\tau\over 2} \right)\right|^2   \sum_{k=0}^m {(-m)_k \over k!}  $$
 
$$\times {\left( \gamma+1-{i\tau\over 2} \right)_k \left( \gamma+1-{i\tau\over 2} \right)_k\over (1+\alpha)_k}\  W_{ - (1/2+\gamma+k), i\tau/2} \left({\eta^2\over 4(1-\mu)}\right).\eqno(2.19)$$ 
The limit case $\mu=1$ is treated, employing the equality  (cf. [12])

$$ \int_0^\infty x^{s-1}  K_{i\tau}\left(\eta\sqrt x\right) dx = \eta^{-2s} 2^{2s-1} \Gamma\left( s + {i\tau\over 2}\right) \Gamma\left( s - {i\tau\over 2}\right),\ {\rm Re}(s) > 0.\eqno(2.20)$$
Hence, correspondingly,

$$G_m\left(\tau; \alpha,\gamma,0, \eta\right) =  {1\over 2}  \sum_{k=0}^m {(-1)^k \over k!} \binom{m+\alpha}{m-k}  \left({4\over \mu^2} \right)^{\gamma+k+1}  \Gamma\left( k+\gamma+1 + {i\tau\over 2}\right) $$

$$\times \Gamma\left( k+\gamma+1 - {i\tau\over 2}\right)=  {(1+\alpha)_m \over 2 m!} \left({4\over \eta^2} \right)^{\gamma+1} \left|\Gamma\left( \gamma+1 + {i\tau\over 2}\right)\right|^2 $$

$$\times  {}_3F_1 \left(  -m,\  \gamma+1+{i\tau\over 2},\ \gamma+1- {i\tau\over 2} ; \  1+ \alpha;\  {4\over \eta^2} \right).\eqno(2.21)$$

Now, considering the Continuous Dual Hahn orthogonality, we write it in the form (see (1.30))

$${1\over \pi^2} \int_0^\infty \tau\sinh(\pi\tau) \left|\Gamma\left(c+{i\tau\over 2}\right)\right|^2  (F_{2,1/2}\  f_n)\left( {\tau^2\over 4}\right) (F_{2,1/2}\  g_m)\left({\tau^2\over 4}\right) d\tau $$

$$= 4 \delta_{n,m} \  n! \Gamma(n+a+b)\Gamma(n+a+c)\Gamma(n+b+c),\eqno(2.22)$$
where via (2.20)

$$f_n(x)=  2 (a+b)_n(a+c)_n\  x^{a-1}  {}_1F_2 \left( -n;  a+b, a+c; x\right),\eqno(2.23)$$

$$g_m(x)=  2 (a+b)_m(b+c)_m  x^{b-1}  {}_1F_2 \left( -m;  a+b, b+c; x\right).\eqno(2.24)$$
We will show that (2.22) generates the convolution orthogonality of functions $f_n,\ g_m$ related to the modified Kontorovich-Lebedev transform (1.30) for $\alpha=2, \beta= 1/2$. To do this,  we will modify the convolution (1.14). Indeed, making use simple substitutions,  we define the modified convolution as follows

$$(f\hat{*}g)(x)= {1\over 4x} \int_0^\infty\int_0^\infty e^{- (y+t) \sqrt {{x\over yt}}  - \sqrt{{yt\over x}} }f(y)g(t) dydt,\quad x >0\eqno(2.25)$$
and (1.33) implies an analog of the equality (1.23)

$$\int_0^\infty  (f\hat{*}g)(x) \omega(x) dx =  {1\over \pi^2} \int_0^\infty \tau\sinh(\pi\tau) \hat{q}(\tau)  (F_{2,1/2}  f)\left(\tau\right) (F_{2,1/2} g)\left(\tau\right) d\tau,\eqno(2.26)$$
where

$$\hat{q}(\tau) =  \int_0^\infty   K_{i\tau} (2\sqrt x) \omega(x) {dx\over x}.\eqno(2.27)$$
Hence functions $f_n, \ g_m$ are  orthogonal with respect to convolution (2.25),  and (2.20), (2.22), (2.26)  yield the equality 

$$ \int_0^\infty  (f_n\hat{*}g_m)(x) x^c dx = 2 \delta_{n,m} \  n! \Gamma(n+a+b)\Gamma(n+a+c)\Gamma(n+b+c).\eqno(2.28)$$
In order to employ convolution (1.14) and Parseval-type equality (1.23) we will appeal to (1.1), (2.7) and Entry 3.14.3.2 in [3] $(\mu >0)$

$$\int_0^\infty  x^{s-1} e^{\mu x} K_{i\tau} (\mu x) dx  = {\cosh(\pi\tau) \over \sqrt \pi (2\mu)^s} \Gamma(s-i\tau)\Gamma(s+i\tau) \Gamma\left({1\over 2}- s\right),\  0 < {\rm Re} s < {1\over 2}.\eqno(2.29)$$
Then the Continuous Dual Hahn orthogonality can be written as follows

$${2\over \pi^2} \int_0^\infty \tau\sinh(\pi\tau) q(\tau)  (F f_n)(\tau)(F g_m)(\tau)  d\tau $$

$$=   {n!\over 2^c \sqrt\pi } \frac{\Gamma(n+a+b)\Gamma(n+a+c)\Gamma(n+b+c)\Gamma(1/2-c)} {\Gamma(a+1/2) \Gamma(b+1/2)}\ \delta_{n,m},\eqno(2.30)$$
where  $a,b, > 0, \ 0 < c < 1/2$ and 

$$f_n(x)= {2^a \over \sqrt\pi} \ x^{a-1} e^{-x}  (a+b)_n(a+c)_n \  {}_2F_2 \left( -n, a+ 1/2; a+b, a+c; 2x \right),\eqno(2.31)$$

$$g_m(x)= {2^b \over \sqrt\pi} \ x^{b-1} e^{-x}  (a+b)_m(b+c)_m\  {}_2F_2 \left( -m, b+ 1/2; a+b, b+c; 2x \right),\eqno(2.32)$$

$$q(\tau)= {\cosh(\pi\tau) \over 2^c \sqrt\pi} \left| \Gamma(c+i\tau)\right|^2  \Gamma\left({1\over 2}-c\right)=  \int_0^\infty  x^{c-1} e^{ x} K_{i\tau} ( x) dx.\eqno(2.33)$$
Hence since $f_n, g_m$ given by (2.31), (2.32) evidently satisfy conditions of Theorem 1, the corresponding Parseval-type equality (1.17) for the convolution $f_n*g_m$ holds. Multiplying both sides by $e^x x^c$ and integrating over $\mathbb{R}_+$, we obtain

$$  \int_0^\infty (f_n*g_m)(x) e^x x^c  dx = {2\over \pi^2} \int_0^\infty \tau\sinh(\pi\tau) q(\tau)  (F f_n)(\tau)(F g_m)(\tau)  d\tau, $$
where $q$ is defined by (2.33) and the interchange of the order of integration is ensured by Fubini's theorem owing to the estimate (see (1.12))

$$  \int_0^\infty \tau\sinh(\pi\tau) \left|  (F f_n)(\tau)(F g_m)(\tau) \right|  \int_0^\infty  x^{c-1} e^{ x} \left|K_{i\tau} ( x)\right| dx d\tau$$ 

$$\le  \int_0^\infty  x^{c-1} e^{ x} K_{0} ( x) dx  \int_0^\infty \tau\sinh(\pi\tau) \left|  (F f_n)(\tau)(F g_m)(\tau) \right| d\tau$$

$$\le  { \pi\sqrt\pi \left| \Gamma(c)\right|^2\over 2^{c+1} } ||f_n||_{L_2\left(\mathbb{R}_+: xdx\right)} \  ||g_m||_{L_2\left(\mathbb{R}_+: xdx\right)} < \infty.$$
Consequently, combining with (2.30), we end up with  the following convolution orthogonality 

$$  \int_0^\infty (f_n*g_m)(x) e^x x^c  dx = {n!\over 2^c \sqrt\pi } \frac{\Gamma(n+a+b)\Gamma(n+a+c)\Gamma(n+b+c)\Gamma(1/2-c)} {\Gamma(a+1/2) \Gamma(b+1/2)}\ \delta_{n,m}.\eqno(2.34)$$
Let us consider the Wilson orthogonality and Wilson's polynomials (2.2). We will show that it generates the convolution orthogonality 
by means of the modified Kontorovich-Lebedev transform (1.30) for $\alpha=2, \beta= 1/2$. In fact, using  Entry  3.14.18.4 in [3]

$$\int_0^\infty K_{i\tau} (2\sqrt x) K_\nu(2\sqrt x) x^{s-1} dx=   {1\over 4 \Gamma(2s)}   \Gamma\left( s+ {\nu+i\tau\over 2} \right)
\Gamma\left( s- {\nu+i\tau\over 2} \right)$$

$$\times  \Gamma\left( s+ {\nu- i\tau\over 2} \right)\Gamma\left( s- {\nu-i\tau\over 2} \right),\eqno(2.35)$$
we derive for the case of Wilson's orthogonality similar to (2.22)

$${1\over \pi^2} \int_0^\infty  \tau \sinh(\pi\tau) \left|  \Gamma\left(c+{i\tau\over 2}\right)\Gamma\left(d+{i\tau\over 2}\right) \right|^2 (F_{2,1/2}\  f_n)\left( {\tau^2\over 4}\right) (F_{2,1/2}\  g_m)\left({\tau^2\over 4}\right) d\tau $$

$$= 4 \delta_{n,m} \frac{ n! \Gamma(n+a+b)\Gamma(n+a+c)\Gamma(n+a+d)\Gamma(n+b+c)\Gamma(n+b+d)\Gamma(n+c+d)}
{\Gamma(2n+a+b+c+d) (n+a+b+c+d-1)_n},\eqno(2.36)$$
where, recalling (2.20),   

$$f_n(x)= 2 x^{a-1}  (a+b)_n(a+c)_n (a+d)_n$$

$$\times  {}_2F_3 \left( -n, n+a+b+c+d-1; a+b, a+c, a+d; x\right),\eqno(2.37)$$

$$g_m(x)= 2 x^{b-1}  (a+b)_n(b+c)_n (b+d)_n$$

$$\times  {}_2F_3 \left( -n, n+a+b+c+d-1; a+b, b+c, b+d; x\right).\eqno(2.38)$$
Therefore,  appealing to (2.26), we find that $f_n, g_m$ given by (2.37), (2.38) are orthogonal with respect to convolution (2.25),  and (2.35), (2.36) imply the equality 

$$\int_0^\infty  (f\hat{*}g)(x) K_{c-d}\left(2\sqrt x\right) x^{(c+d)/2}  dx $$

$$=   \delta_{n,m} \frac{ n! \Gamma(n+a+b)\Gamma(n+a+c)\Gamma(n+a+d)\Gamma(n+b+c)\Gamma(n+b+d)(c+d)_n}{ \Gamma(2n+a+b+c+d) (n+a+b+c+d-1)_n}.\eqno(2.39)$$
On the other hand, doing analogously to (2.30), we involve convolution (1.14) and equalities (1.16), (2.7), (2.29), writing the right-hand side of (1.23) accordingly 

$$ {2\over   \pi^2} \int_0^\infty \tau\sinh(\pi\tau)  q(\tau)  (Ff_n)(\tau) (Fg_m)(\tau) d\tau $$

$$=  {\delta_{n,m} \over 2^{c+d} } \frac{ n! \Gamma(n+a+b)\Gamma(n+a+c)\Gamma(n+a+d)\Gamma(n+b+c)\Gamma(n+b+d)\Gamma(n+c+d)\Gamma(1/2-c)}
{\Gamma(1/2+d) \Gamma(2n+a+b+c+d) (n+a+b+c+d-1)_n},\eqno(2.40)$$
where $0 < c < 1/2$ and 

$$q(\tau)= \int_0^\infty K_{i\tau} (x) \left( t^{c-1} e^t * t^{d-1} e^{-t} \right)(x)  dx,\eqno(2.41)$$

$$f_n(x)=  {2^a\over \sqrt \pi}  x^{a-1} e^{-x} (a+b)_n(a+c)_n (a+d)_n$$

$$\times  {}_3F_3 \left( -n, n+a+b+c+d-1, a+ 1/2; a+b, a+c, a+d; 2x\right),\eqno(2.42)$$

$$g_m(x)=    {2^b\over \sqrt \pi}  x^{b-1} e^{-x} (a+b)_m(b+c)_m (b+d)_m$$

$$\times  {}_3F_3 \left( -m, m+a+b+c+d-1, b+ 1/2; a+b, b+c, b+d; 2x\right).\eqno(2.43)$$
The convolution  $\left( t^{c-1} e^t * t^{d-1} e^{-t} \right)(x)$ in (2.41) can be treated, taking relation  2.2.1.8  in [3] 

$$\int_0^\infty x^{s-1} e^{-ax- b/x} dx = 2 \left({b\over a}\right)^{s/2} K_s\left(2\sqrt{ab}\right).\eqno(2.44)$$
Then we get from (1.14)

$$ \left( t^{c-1} e^t * t^{d-1} e^{-t} \right)(x) = x^{d-1}  \int_0^\infty  {K_d(x+y)\over (x+y)^d}\ e^y   y^{c+d-1}  dy.\eqno(2.45)$$
Consequently, we arrive at the convolution orthogonality of functions (2.42), (2.43), obtaining from (2.40) and (1.23) the equality 

$$\int_0^\infty (f_n*g_m)(x) \omega(x) dx $$

$$=  {\delta_{n,m} \over 2^{c+d} } \frac{ n! \Gamma(n+a+b)\Gamma(n+a+c)\Gamma(n+a+d)\Gamma(n+b+c)\Gamma(n+b+d)\Gamma(n+c+d)\Gamma(1/2-c)} {\Gamma(1/2+d) \Gamma(2n+a+b+c+d) (n+a+b+c+d-1)_n},\eqno(2.46)$$
where

$$\omega(x)= x^{d}  \int_0^\infty  {K_d(x+y)\over (x+y)^d}\ e^y   y^{c+d-1}  dy,\quad x >0.$$

\section{Orthogonal and $d$-orthogonal polynomials of the Prudnikov type}

In this section we shall show how the Prudnikov-type orthogonal  and $d$-orthogonal polynomials (see [4], [14], [15]) , which are associated with the scaled Macdonald functions $\rho_\nu(x)= 2 x^{\nu/2} K_\nu(2\sqrt x)$,  generate new orthogonal systems of functions by virtue of the modified Kontorovich-Lebedev transform (1.30) $F_{2,1/2} f$.  Precisely,  let us consider the following orthogonalities

$$\int_0^\infty  P_n^{\nu,\alpha} (x)P_m^{\nu,\alpha} (x)  x^\alpha \rho_\nu(x) dx=  \delta_{n,m} , \quad \nu \ge 0,\ \alpha > 0,\eqno(3.1)$$

$$\int_{0}^{\infty}  Q^\nu_{m} (x) Q^\nu_{n}(x) e^{-x} \rho_\nu(x) dx = \delta_{n,m}, \quad \nu > 0,\eqno(3.2)$$

$$\int_{0}^{\infty}  q^\nu_{m} (x) q^\nu_{n}(x) e^{- 1/x} \rho_\nu(x) {dx\over x} = \delta_{n,m}, \quad \nu > 0,\eqno(3.3)$$
and $d$-orthogonality  conditions $(\alpha >0)$

$$\int_0^\infty  p_{2n}^{\nu,\alpha} (x)  \rho_\nu(x)  x^{\alpha+m} dx= 0,\quad m= 0,1, 2, \dots, n-1,\eqno(3.4)$$

$$\int_0^\infty  p_{2n}^{\nu,\alpha} (x)  \rho_{\nu+1} (x)  x^{\alpha+m} dx= 0,\quad m= 0,1, 2, \dots, n-1,\eqno(3.5)$$

$$\int_0^\infty  p_{2n+1}^{\nu,\alpha} (x) \rho_\nu(x)  x^{\alpha+m} dx= 0,\quad m= 0,1, 2, \dots, n,\eqno(3.6)$$

$$\int_0^\infty  p_{2n+1}^{\nu,\alpha} (x)  \rho_{\nu+1} (x)  x^{\alpha+m} dx= 0,\quad m= 0,1, 2, \dots, n-1.\eqno(3.7)$$
Then formulas (2.20),  (2.35), Parseval equality (1.33)  and (3.1) generate the following orthogonality 

$$ \int_0^\infty   S_n^{\nu,\alpha} (2\tau) U_m^{\nu,\alpha} (2\tau) \left| \frac{ \tau \Gamma\left(\nu+ \alpha+i\tau\right) \Gamma\left(\alpha+ i\tau\right)}{\Gamma\left(i\tau+1/ 2\right)}\right|^2  d\tau $$

$$= {1\over 2} \  \Gamma(2\alpha+\nu) \delta_{n,m},\eqno(3.8)$$
where polynomials  $S_n^{\nu,\alpha} (2\tau),\  U_m^{\nu,\alpha} (2\tau)$ are defined accordingly

$$ S_n^{\nu,\alpha} (2\tau) = \sum_{k=0}^n a_{n,k} (1+i\tau)_k (1-i\tau)_k,\eqno(3.9)$$

$$U_m^{\nu,\alpha} (2\tau) = \sum_{k=0}^m {a_{m,k}\over (2\alpha+\nu)_{2k} } (\alpha+\nu+ i\tau)_k (\alpha +\nu-i\tau)_k (\alpha+ i\tau)_k (\alpha -i\tau)_k,\eqno(3.10)$$
and coefficients $a_{n,k}$ of the Prudnikov polynomials $P_n^{\nu,\alpha}$ are given explicitly in [14].  Concerning the Prudnikov-type polynomials $Q_n^\nu$, we appeal again to (1.33),  (2.35) and Entry 3.14.3.10 in [3] (cf. (2.18)) to get from (3.2) the equality 

$$  \int_0^\infty \tau\sinh(\pi\tau) F_{2,1/2} \left( Q^\nu_{n}(x) e^{-x} x^{\nu/2-1}\right)(\tau) $$

$$\times   F_{2,1/2} \left( Q^\nu_{m} (x)  K_\nu(2\sqrt x)  \right) (\tau) d\tau= {\pi^2\over 2} \delta_{n,m},\eqno(3.11) $$
where 

$$F_{2,1/2} \left( Q^\nu_{n}(x) e^{-x}   x^{\nu/2-1} \right)(\tau) = {\sqrt e\over 2} \left|\Gamma\left( {i\tau+\nu\over 2} \right)\right|^2 \sum_{k=0}^n b_{n,k} \left( {\nu+ i\tau\over 2} \right)_k \left( {\nu- i\tau\over 2} \right)_k $$

$$\times W_{(1-\nu)/2-k, i\tau/2} (1),\eqno(3.12)$$

$$F_{2,1/2} \left( Q^\nu_{m} (x)  K_\nu(2\sqrt x) \right) (\tau) =  {1\over 4 } \left|\Gamma\left(1+ {\nu+ i\tau\over 2} \right) \Gamma\left( 1+{\nu-i\tau\over 2} \right)\right|^2$$

$$\times \sum_{k=0}^m {b_{m,k} \over (2)_{2k}}  \left(1+{\nu+ i\tau\over 2} \right)_k \left(1-{\nu+i\tau\over 2} \right)_k \left(1+{\nu-i\tau\over 2} \right)_k \left(1- {\nu-i\tau\over 2} \right)_k,\eqno(3.13)$$
and coefficients $b_{m,k}$ of the polynomials $Q_n^\nu$ are calculated explicitly in [15].   Analogously we treat the orthogonality (3.3).  However, to do this,  we will need the formula (see [3], Entry 3.14.3.13 

$$2 \int_0^\infty e^{-1/x} K_{i\tau} (2\sqrt x) x^{s-1} dx =  \Gamma\left(s+ {i\tau\over 2} \right) \Gamma\left(s- {i\tau\over 2} \right)
{}_0F_2 \left( 1-s+ {i\tau\over 2},\  1-s- {i\tau\over 2};\ -1\right) $$

$$+   \Gamma\left(- s- {i\tau\over 2} \right) \Gamma\left(- i\tau \right) {}_0F_2 \left( 1+s+ {i\tau\over 2},\  1+ i\tau;\ -1\right) $$

$$+  \Gamma\left(- s + {i\tau\over 2} \right) \Gamma\left( i\tau \right) {}_0F_2 \left( 1+s- {i\tau\over 2},\  1- i\tau;\ -1\right) .$$
Hence, recalling (1.33) and (2.35), we derive 

$$ \int_0^\infty \tau\sinh(\pi\tau) F_{2,1/2} \left( q^\nu_{n}(x) e^{-1/x} x^{-2} \right)(\tau)   F_{2,1/2} \left( q^\nu_{m} (x) \rho_\nu(x) \right) (\tau) d\tau= \pi^2  \delta_{n,m},\eqno(3.14) $$
where

$$ F_{2,1/2} \left( q^\nu_{n}(x) e^{-1/x} x^{-2} \right)(\tau) = {1\over 2} \left|\Gamma\left( {i\tau\over 2} -1\right)\right|^2  \sum_{k=0}^n c_{n,k} \left( {i\tau\over 2}-1 \right)_k \left(- {i\tau\over 2} -1\right)_k $$

$$\times {}_0F_2 \left(2 -k+ {i\tau\over 2},\  2-k- {i\tau\over 2};\ -1\right) + {\Gamma\left(-  i\tau\right) \over 2} \sum_{k=0}^n c_{n,k}  \Gamma\left(1- k- {i\tau\over 2} \right) $$

$$\times {}_0F_2 \left( k+ {i\tau\over 2},\  1+ {i\tau\over 2};\ -1\right) + {\Gamma\left(  i\tau\right) \over 2} \sum_{k=0}^n c_{n,k}  \Gamma\left(1- k+{i\tau\over 2} \right) $$

$$\times {}_0F_2 \left( k- {i\tau\over 2},\  1- {i\tau\over 2};\ -1\right),\eqno(3.15)$$

$$F_{2,1/2} \left( q^\nu_{m} (x) \rho_\nu(x) \right) (\tau) = {1\over 2 \Gamma(2+\nu) } \left|\Gamma\left(1+\nu+ {i\tau\over 2} \right) \Gamma\left( 1+{i\tau\over 2} \right)\right|^2$$

$$\times \sum_{k=0}^m {c_{m,k} \over (2+\nu)_{2k}}  \left(1+\nu+ {i\tau\over 2} \right)_k \left(1+\nu- {i\tau\over 2} \right)_k \left(1+{i\tau\over 2} \right)_k \left(1- {i\tau\over 2} \right)_k,\eqno(3.16)$$
and coefficients $c_{m,k}$ of the polynomials $q_n^\nu$ are calculated explicitly in [15].  Finally, we generate new system of $d$-orthogonal polynomials, exploring the orthogonality conditions  (3.4)-(3.7).  We will proceed this with the use of the explicit formula for polynomials $p_n$ from [4]

$$p_n^{\nu,\alpha} (x) = (-1)^n (1+\alpha)_n(1+\alpha+\nu)_k \  {}_1F_2 \left( -n;\ 1+\alpha,\ 1+\alpha+\nu;\  x\right).$$
Then we derive from (3.4)-(3.7), correspondingly, employing (1.33), (2.20), (2.35),  

$$\int_0^\infty  V_{2n}^{\nu,\alpha} (\tau^4 )  \left|{\tau \Gamma\left(1+\nu+ i\tau \right)\Gamma\left(\alpha+ i\tau +m\right)\over \Gamma(i\tau+1/2) }\right|^2 d\tau = 0,\ m= 0,1, 2, \dots, n-1,\eqno(3.17)$$

$$\int_0^\infty  S_{2n}^{\nu,\alpha} (\tau^4 )  \left|{\tau \Gamma\left(2+\nu+ i\tau \right)\Gamma\left(\alpha+ i\tau +m\right)\over \Gamma(i\tau+1/2) }\right|^2 d\tau = 0,\  m= 0,1, 2, \dots, n-1,\eqno(3.18)$$

$$\int_0^\infty  V_{2n+1}^{\nu,\alpha} (\tau^4 )  \left|{\tau \Gamma\left(1+\nu+ i\tau \right)\Gamma\left(\alpha+ i\tau +m\right)\over \Gamma(i\tau+1/2) }\right|^2 d\tau = 0,\ m= 0,1, 2, \dots, n,\eqno(3.19)$$

$$\int_0^\infty  S_{2n+1}^{\nu,\alpha} (\tau^4 )  \left|{\tau \Gamma\left(2+\nu+ i\tau \right)\Gamma\left(\alpha+ i\tau +m\right)\over \Gamma(i\tau+1/2) }\right|^2 d\tau = 0,\ m= 0,1, 2, \dots, n-1,\eqno(3.20)$$
where 

$$V_n^{\nu,\alpha} (\tau^4)=  (-1)^n (1+\alpha)_n(1+\alpha+\nu)_k $$

$$\times   {}_5F_4 \left( -n, 1+\nu+ i\tau, 1+\nu- i\tau,  1+ i\tau,  1- i\tau ; {\nu+2\over 2},  {\nu+3\over 2}, \ 1+\alpha,  1+\alpha+\nu;\  {1\over 4}\right), $$

$$S_n^{\nu,\alpha} (\tau^4)=  (-1)^n (1+\alpha)_n(1+\alpha+\nu)_k $$

$$\times   {}_5F_4 \left( -n, 2+\nu+ i\tau, 2+\nu- i\tau,  1+ i\tau,  1- i\tau ; {\nu+3\over 2},  {\nu+4\over 2}, \ 1+\alpha,  1+\alpha+\nu;\  {1\over 4}\right).$$

\bigskip
\centerline{{\bf Acknowledgments}}
\bigskip
The work was partially supported by CMUP (UID/MAT/00144/2019),  which is funded by FCT (Portugal) with national (MEC),   European structural funds through the programs FEDER  under the partnership agreement PT2020. 

\bigskip
\centerline{{\bf References}}
\bigskip
\baselineskip=12pt
\medskip
\begin{enumerate}

\item[{\bf 1.}\ ]
 W.A. Al-Salam, M. Ismail, Polynomials orthogonal with respect to discrete convolution, {\it  J. Math. Anal. Appl.,}  (1976), {\bf 55},  125-139. 

\item[{\bf 2.}\ ]
 V. Amerbaev, Z. Naurzbaev, Polynomials orthogonal with respect to  convolution, {\it  Izv. Akad. Nauk Kazah. SSR Ser. Fiz.-Mat. Nauk,}  (1965), {\bf 3},   70-78 (in Russian).

\item[{\bf 3.}\ ]  Yu.A. Brychkov, O.I. Marichev, N.V. Savischenko,   {\it Handbook of Mellin Transforms.}  
Advances in Applied Mathematics, CRC Press,  Boca Raton, 2019.

\item[{\bf 4.}\ ]  E. Coussement, W. Van Assche, Some properties of multiple orthogonal polynomials associated with Macdonald functions,     J. Comput. Appl. Math. 133 (2001),  253-261.

\item[{\bf 5.}\ ] A. Erd\'elyi, W. Magnus, F.  Oberhettinger, and F.G. Tricomi,   {\it Higher Transcendental Functions.} Vols. I and
II, McGraw-Hill, New York, London and Toronto, 1953.

\item[{\bf 6.}\ ]
 T.H. Koornwinder,  Special orthogonal polynomial systems mapped onto each other by the Fourier-Jacobi transform,
 {\it  Lecture Notes in Math.,} {\bf 1171}, 174-183, Springer, Berlin, 1985.

 \item[{\bf 7.}\ ]
 A.F Loureiro, S. Yakubovich, The Kontorovich-Lebedev transform as a map between $d$-orthogonal polynomials, {\it Stud. Appl. Math.,}  (2013), {\bf 131}, 229-265. 

\item[{\bf 8.}\ ]
 A.F Loureiro, S. Yakubovich, Central factorials under the Kontorovich-Lebedev transform of polynomials, {\it  Integral Transforms Spec. Funct.,}  (2013), {\bf 24}, 217-238. 

\item[{\bf 9.}\ ]
 M. Masjed-Jamei, W. Koepf,   Two classes of special functions using Fourier transforms of generalized ultraspherical and generalized Hermite polynomials, {\it  Proc. Amer. Math. Soc.,}  (2012), {\bf 140},  N 6, 2053-2063. 

\item[{\bf 10.}\ ]
 E.C. Titchmarsh,  {\it An Introduction to the Theory of Fourier Integrals}, Clarendon Press, Oxford, 1937.

\item[{\bf 11.}\ ]
 J.A. Wilson, Some hypergeometric orthogonal polynomials, {\it  SIAM J. Math. Anal.,}  (1980), {\bf 11}, N 4, 690-701.

\item[{\bf 12.}\ ]  S. Yakubovich, Yu. Luchko, {\it The hypergeometric approach to integral transforms and convolutions. Mathematics and its Applications, 287,}  {\it Kluwer Academic Publishers Group}, Dordrecht, 1994.

\item[{\bf 13.}\ ] S. Yakubovich, {\it Index Transforms}, World Scientific Publishing Company, Singapore, New Jersey, London and
Hong Kong, 1996.

\item[{\bf 14.}\ ] S. Yakubovich,  Orthogonal polynomials with ultra-exponential weight functions: an explicit solution to the Ditkin-Prudnikov problem.  ArXiv:1811.03475.

\item[{\bf 15.}\ ] S. Yakubovich, Orthogonal polynomials with the Prudnikov-type weights.  ArXiv:1902.06227.

\end{enumerate}

\vspace{5mm}

\noindent S.B.Yakubovich\\
Department of  Mathematics,\\
Faculty of Sciences,\\
University of Porto,\\
Campo Alegre st., 687\\
4169-007 Porto\\
Portugal\\
E-Mail: syakubov@fc.up.pt\\

\end{document}